\documentstyle[12pt]{article}
\begin{document}
\title{Noncommutative Vieta Theorem and Symmetric Functions}
\author{Israel Gelfand \and Vladimir Retakh\\
Department of Mathematics, Rutgers University\\
New Brunswick, NJ 08903}
\date{July 10, 1995}
\maketitle

{\it This paper is submitted to Gelfand Seminar Volume 1994-95}
\vskip .5cm
There are two ways to generalize basic constructions of commutative algebra
for a noncommutative case. More traditional way is to define commutative
functions like trace or determinant over noncommuting variables. Beginning
with [6] this approach was widely used by different
authors, see for example [5], [15], [14], [12], [11], [7].

However, there is another possibility to work with purely noncommutative
objects without using trace or determinant or passing to a quotient space or
quotient algebra started in [9] and [10]. Let us compare these two approaches
on a simplest example -
a classical Vieta theorem which, of course, is a starting point in a theory of
symmetric functions.

Consider an algebraic equation
\begin{equation}
x^n+a_1x^{n-1}+a_2x^{n-2}+\ldots +a_n=0
\end{equation}

In a noncommutative case first formulas expressing coefficients
$a_1, a_2, \ldots , a_n$ via solutions $x_1, x_2, \ldots , x_n$
of the equation (1) have appeared in [8], Section 7.1 .
This paper heavily used a theory of quasideterminants over
(non)commutative rings developed in [9], [10]. The expressions for
$a_1, \ldots , a_n$ were given via ratios of Vandermonde quasideterminants
depending of variables $x_1, \ldots , x_n$. These expressions in general
are rational functions of $x_1,\ldots , x_n$. One must use then nontrivial
determinant identities to obtain classical Vieta formulas in a
commutative case.

In an interesting paper [7] Fuchs and Schwarz tried to give an analogue of
classical
Vieta formulas in a noncommutative case. Let the coefficients of the
equation (1) belong to an algebra $R$ over a field $k$ and $x_1, \ldots ,
x_n$ be a set of independent
solutions of (1) (For a matrix case it means that the corresponding
bloc Vandermonde determinant is not equal to
zero). Fuchs and Schwarz proved

{\bf Theorem 1}
{\it If there is an additive morphism}
${tr}:R\rightarrow k$
{\it satisfying the condition}
${tr}\ uv={tr}\ vu$ {\it for any}
$u,v\in F$,
{\it then}
$$
  {tr}\ a_1=-({tr}\ x_1+\ldots +{tr}\ x_n).
$$
{\it If there is a multiplicative morphism}
${det}:R\rightarrow k$,
{\it then}
$$
{det}\ a_n={det}(-x_1)\cdot \ldots \cdot {det}\left(
-x_n\right) .
$$

This result was proved in [7] when $a_1, \ldots , a_n$ and $x_1,\ldots ,x_n$
are just complex
matrices. Then the authors used the Amitsur theorem that
identities which hold
for complex matrices hold also in arbitary associative rings with units.
There are no similar formulas in [7] for intermediate coefficients
$a_2, \ldots , a_{n-1}$.

In this paper we will give much more general version of noncommutative
Vieta theorem. It does not require the existence of trace or determinant and
also give formulas for intermiediate coefficients.

Namely, for a ``generic'' set of solutions $x_1,\ldots , x_n$ of the equation
(1) over a (noncommutative) skew-field we will construct a set of rational
functions $v_1,\ldots , v_n$ depending of $x_1,\ldots , x_n$ and a set
of variables

\begin{equation}
y_k=v_kx_kv_k^{-1},
\end{equation}
where $k=1,2,\dots ,n$.
We call $v_i$'s Vandermonde quasideterminants (see section 2). Our first
main result is  that

\[
a_k=(-1)^k\sum_{1\leq i_1<i_2\ldots <i_k\leq n}y_{i_k}\cdot \ldots \cdot
y_{i_1}
\]
for $k=1,2,...,n$.

In particular,
\begin{eqnarray*}
a_1 &=&-(y_1+y_2+\ldots +y_n), \\
a_2 &=&\sum_{i<j}y_jy_i ,\\
a_n &=&(-1)^ny_ny_{n-1}\cdot \ldots \cdot y_1.
\end{eqnarray*}

Theorem 1 immediately follows from our statement and formulas (2). We do not
use Amitsur Theorem here. Our
proof is based on "honest" algebraic computations using
quasideterminant identities. For these reasons
we are working inside free skew-fields generated by a finite set of
noncommutative variables over a commutative field.

The expressions $a_k$ or $\Lambda _k=(-1)^ka_k$ for $k=1,\dots ,n$
are symmetric functions of $x_1,\dots ,x_n$. Following a general line
of [8] consider a free associative algebra ${\bf Symm}$
over a noncommutative field of characteristics zero generated by
$\Lambda _1, \Lambda _2,\dots , \Lambda _n$. Each element of this algebra
may be viewed as a polynomial of $y_1,\dots ,y_n$ as well as a
rational function of $x_1,\dots ,x_n$.

{\bf Theorem 5}. {\it A polynomial} $P$ {\it of} $y_1,\dots ,y_n$
{\it is symmetric in} $x_1,\dots ,x_n$ {\it if and only if} $P$ {\it
belongs to the algebra} ${\bf Symm}$.

In other words, this theorem shows that a ``real'' symmetric functions
$P$ is an ``abstract'' symmetric function in a sense of [8].
According to Section 7.3 of [8], we describe a basis in these
functions.

Let $w = y_{i_1}\dots y_{i_m}$ be a word. An integer $k$ is called
a {\it descent} of $w$ if $ 1\leq k\leq m-1$ and $i_k$
is greater than
$i_{k+1}$.

For any set $J=(j_1,\dots ,j_k)$ of nonnegative integers consider a
function
$$R_J(y_1,\dots ,y_n) = \sum  y_{i_1}\dots y_{i_m}, $$

where $m = j_1+\dots +j_k$ and the sum is running over all words
$w= y_{i_1}\dots y_{i_m}$ where descents are precisely $j_1, j_1+j_2,
\dots , j_1+\dots +j_{k-1}$.

Such functions were called {\it ribbon Schur functions} in [8].
{}From Proposition 7.15 in [8], Section 7.3 it follows that functions
$R_J$ form a linear base in ${\bf Symm}$. It means that the base in
${\bf Symm}$ is parametrized by sequences of nonnegative integers.
We recall that in a commutative case the well-known base of
classical Schur functions is parametrized by {\it weakly increasing}
sequences of nonnegative integers.

Note also that functions
$$ S_k(y_1,\dots ,y_n)=\sum y_{i_1}y_{i_2}\dots y_{i_k}$$
where the sum is running over all $i_1\leq i_2 \leq \dots \leq i_k$
constitue a special form of functions $R_J$. Functions $S_k$ are
analogues of complete symmetric functions in a commutative case. They were
considered in [8], Section 7.3.

{\bf 1}. We recall here a notion of quasideterminant defined in [9], [10].
Let $A=\left\{ \left\| a_{ij}\right\| , i\in I, j\in J\right\} $ be a square
matrix of order $n=\left| I\right| =\left| J\right| $ with formal
noncommutative entries $a_{ij}$. For $p\in I,$ $q\in J$ denote by $A^{pq}$
the submatrix $\left\{ \left\| a_{ij}\right\| , i\in I-\{p\}, j\in
J-\{q\}\right\} $ of $A$. Let $F$ be a free skew-field defined by formal
variables $a_{ij}$.

{\bf Definition}
{\it The formula}
\[
\left| A\right| _{pq}=a_{pq}-\sum_{i\in I-\{p\}, j\in J-\{q\}}
a_{pj}\left| A^{pq}\right| _{ij}^{-1}a_{iq}
\]
({\it which reduces to $\left| A\right| _{pq}=a_{pq}$ if $n=1$) defines
inductively $n^2$ quasideterminants $\left| A\right| _{pq}$ of the matrix $A$.}

This definition is also valid for a generic matrices over a skew-field, i.e.
matrices for which all expressions in the formula for $|A|_{pq}$ are defined.

{\bf Remark}. Free skew-fields were introduced by Amitsur [1] and studied by
Bergman [2]
and P. Cohn [3],  [4] who characterized them as universal skew-field of
fractions of the ring of noncommutative polynomials; they are universal
objects in the cathegory whose morphisms are specializations. A reader who
is unfamiliar with this subject may just consider our expressions in a generic
case.

{\bf Example}
For $n=2$ one has four quasideterminants
\begin{eqnarray*}
\left| A\right| _{11} =a_{11}-a_{12}a_{22}^{-1}a_{21},\ \ \left| A\right|
_{12}=a_{12}-a_{11}a_{21}^{-1}a_{22}, \\
\left| A\right| _{21} =a_{21}-a_{22}a_{12}^{-1}a_{11},\ \ \left| A\right|
_{22}=a_{22}-a_{21}a_{11}^{-1}a_{12}.
\end{eqnarray*}

In the commutative case $\left| A\right| _{pq}=\pm {det}A/
{det}A^{pq}$.

{\bf 2}. Let us construct some expressions which we call Vandermonde
quasideterminants. Suppose that an ordered set
$X=\{x_1<x_2<\ldots <x_n\}$ of solutions of the equation (1) over
a skew-field is given.
Consider for $k=2,3,\ldots ,n$ formal expressions
\begin{equation}
v_k=\left|
\begin{array}{cccc}
x_1^{k-1} & x_2^{k-1} & \cdots & x_k^{k-1} \\
x_1^{k-2} & x_2^{k-2} & \cdots & x_k^{k-2} \\
\cdots & \cdots & \cdots & \cdots \\
1 & 1 & \cdots & 1
\end{array}
\right| _{1k}
\end{equation}
We call these expressions {\it Vandermonde quasideterminants}. We will call
the set of solutions $X$ {\it generic } if all $v_i$'s are defined and
invertible.

{\bf Example}
By the definition of quasideterminants:
$$v_2=x_2  - x_3 , $$

$$v_3=x_3^2 - x_1^2(x_1-x_2)^{-1}x_3 - x_1(1-x_2^{-1}x_1)^{-1}$$

$$-x_2^2(x_2-x_1)^{-1}-x_2^2(1-x_1^{-1}x_2)^{-1}=$$

$$x_3^2-(x_2^2 - x_1^2)(x_2 - x_1)^{-1}x_3 -
(x_2^2-x_1x_2)(1-x_1^{-1}x_2)^{-1}$$.

Set $v_1=1$.

For $k=1,2,\ldots ,n$ define
\begin{equation}
y_k=v_kx_kv_k^{-1}.
\end{equation}

{\bf Example}
$$y_1=x_1 , $$
$$y_2=(x_2-x_1)(x_2(x_2-x_1)^{-1} , $$
$$y_3=\{x_3^2 - (x_2^2-x_1^2)(x_2-x_1)^{-1}x_3 -
       (x_2-x_1)(x_2^{-1}-x_1^{-1})^{-1}\}\cdot x_3\cdot $$
$$\cdot  \{x_3^2 - (x_2^2-x_1^2)(x_2-x_1)^{-1}x_3 -
       (x_2-x_1)(x_2^{-1}-x_1^{-1})^{-1}\}^{-1}$$

Now we formulate our main result.

{\bf Theorem 2}
{\it If $\{x_1,\ldots ,x_n\}$ is an ordered generic set of solutions of the
equation (1) over a skew-field then for $k=1,2,\ldots ,n$
\[
a_k=(-1)^k\sum_{1\leq i_1<i_2\ldots <i_k\leq n}y_{i_k}\cdot \ldots \cdot
y_{i_1}
\]
where $y_k$'s are defined by formula (4)}.

In particular,
\begin{eqnarray*}
a_1 &=&-(y_1+y_2+\ldots +y_n), \\
a_2 &=&\sum_{i<j}y_jy_i ,\\
a_n &=&(-1)^ny_ny_{n-1}\cdot \ldots \cdot y_1.
\end{eqnarray*}

Note that each $y_k$ depends of an ordering of $x_1,\ldots , x_n$ but the
expressions for $a_k$'s do not depend of an ordering.

{\bf 3}. Let us illustrate Theorem 2.
For $n=2$
it is easy to check that $x_1$ and $x_2$ are solutions of the equation
\[
x^2+a_1x+a_2=0
\]
where
\[
a_1=-x_1-(x_2-x_1)x_2(x_2-x_1)^{-1},
\]
\[
a_2=(x_2-x_1)x_2(x_2-x_1)^{-1}x_1,
\]
Note that from the formal identity
\[
x_2(x_2-x_1)^{-1}x_1=(x_1^{-1}-x_2^{-1})^{-1}=x_1(x_2-x_1)^{-1}x_2
\]
it follows that $a_1$and $a_2$ do not depend on the ordering of variables $%
x_1$ and $x_2$.

Note also that  function $y_1y_2$ is not symmetric in $x_1, x_2$ (but
$a_2=y_2y_1$ is!).

It is still possible to check "by hands" that $x_1$, $x_2$, $x_3$ are
solutions of the equation
$$x^3-a_1x^2+a_2x +a_3=0$$

where
$$a_1=-(y_1+y_2+y_3), $$
$$a_2=y_2y_1+y_3y_2+y_3y_1$$
$$a_3=-y_3y_2y_1$$
and $y_1$,$y_2$, $y_3$ are given by formulas (4). One can also check that
$a_1$, $a_2$, $a_3$ are symmetric in $x_1$, $x_2$, $x_3$. However, even in
this case better follow a general proof. Such proof uses Theorem 3 below
and quasideterminant identities from [9] and [10].


{\bf 4}. The following result was essentially obtained in [8], Section 7.1.

{\bf Theorem 3}
{\it Let $x_1,\ldots ,x_n$ be a set of independent solutions of the
equation (1).
Then for} $k=1,...,n$
\[
a_k=-\left|
\begin{array}{cccc}
x_1^n & x_2^n & \cdots & x_n^n \\
\cdots & \cdots & \cdots & \cdots \\
x_1^{n-k+1} & x_2^{n-k+1} & \cdots & x_n^{n-k+1} \\
x_1^{n-k-1} & x_2^{n-k-1} & \cdots & x_n^{n-k-1} \\
\cdots & \cdots & \cdots & \cdots \\
1 & 1 & \cdots & 1
\end{array}
\right| _{1n}\left|
\begin{array}{cccc}
x_1^{n-1} & x_2^{n-1} & \cdots & x_n^{n-1} \\
\cdots & \cdots & \cdots & \cdots \\
x_1^{n-k} & x_2^{n-k} & \cdots & x_n^{n-k} \\
\cdots & \cdots & \cdots & \cdots \\
1 & 1 & \cdots & 1
\end{array}
\right| _{kn}^{-1}.
\]

This theorem demonstrates that coefficients $a_k$ are the ratio of two
quasideterminants of order $n$.

Such expressions with
a change a signs were called in [8] elementary symmetric functions in
$x_1,\ldots ,x_n$. It
was proved in [8] and it also follows from the Theorem 2 that these
functions are really symmetric in $x_1,\ldots ,x_n$.

{\bf Example}
For $n=2$%
\[
a_1=-(x_2^2-x_1^2)(x_2-x_1)^{-1},
\]
\[
a_2=-(x_2^2-x_1x_2)(1-x_1^{-1}x_2)^{-1}=-(x_2-x_1)(x_2^{-1}-x_1^{-1})^{-1}.
\]

Theorem 3 follows from the observation that the equation (1) can be written
in a quasideterminant form
\[
\left|
\begin{array}{cccc}
x^n & x_1^n & \cdots & x_n^n \\
x^{n-1} & x_1^{n-1} & \cdots & x_n^{n-1} \\
\cdots & \cdots & \cdots & \cdots \\
1 & 1 & 1 & 1
\end{array}
\right| _{11}=0.
\]

Transformations $x\rightarrow v_nxv_n^{-1}$ where $v_n$ is a Vandermonde
quasideterminant defined by the formula (2) are usefull for a comparison of
solutions of the left equation (1) and the right equation

\begin{equation}
x^n+x^{n-1}a_1+x^{n-1}a_2+\ldots +a_n=0.
\end{equation}
We give here only the most simple example.

{\bf Theorem 4}
{\it Suppose that $x_1,x_2,\ldots ,x_n$ are independent solutions of the left
equation (1). Then $y_n=v_nx_nv_n^{-1}$ is a solution of the right equation
(5)}.

{\it e-mail}: igelfand@math.rutgers.edu, retakh@dimacs.rutgers.edu
\end{document}